\RequirePackage{tikz} 
\RequirePackage{nicematrix}

\documentclass[pdflatex,sn-mathphys-num]{sn-jnl} 

\theoremstyle{thmstyleone}%
\newtheorem{theorem}{Theorem}
\newtheorem{proposition}[theorem]{Proposition}%

\theoremstyle{thmstyletwo}%

\theoremstyle{thmstylethree}%

\newtheorem*{assumption*}{\assumptionnumber}
\providecommand{\assumptionnumber}{}
\makeatletter

\newcommand{\uu}{u}

\newcommand{\vv}{v}

\newcommand{\mmu}{\mu}

\newcommand{\zz}{z}

\newcommand{\bK}{N}

\newcommand{\bC}{O}

\newcommand{\bB}{Z}
\newcommand{\bM}{M}

\newcommand{\bla}{\lambda}

\newcommand{\TT}{\top}

\usepackage{dutchcal}

\usepackage{bm,breqn} 
\usepackage{multirow,tikz}  
\usepackage{float}
\newcommand\scalemath[2]{\scalebox{#1}{\mbox{\ensuremath{\displaystyle #2}}}}
\usepackage{subcaption}
\usepackage{todonotes}
\usepackage[short]{ optidef }
\usepackage{pgfplots}
\usepackage{pgfplotstable}
\usepackage{subcaption}
\usepackage{amsmath,amsfonts}
\usepackage{algorithm,algpseudocode}
\usepackage{comment}
 \usepackage{tabstackengine}
\stackMath
\definecolor{green}{RGB}{34, 139, 34}

\definecolor{clemson}{RGB}{245,102,0}

\definecolor{FFC20A}{RGB}{255, 194, 10} 
\definecolor{0C7BDC}{RGB}{12, 123, 220} 
\definecolor{E66100}{RGB}{230, 97, 0} 

\begin{document}

\title[A parallel-in-time multigrid preconditioner for optimal control]{A Parallel-in-Time Multigrid Preconditioner for Optimal Control}


\author[1]{\small\fnm{Radoslav} \sur{Vuchkov}}\email{rgvuchk@sandia.gov}

\author[2]{\small\fnm{Eric C.} \sur{Cyr}}\email{eccyr@sandia.gov}

\author[1]{\small\fnm{Aurya} \sur{Javeed}}\email{asjavee@sandia.gov}

\author[1]{\small\fnm{Denis} \sur{Ridzal}}\email{dridzal@sandia.gov}

\affil[1]{\small \orgdiv{Optimization \& Uncertainty Quantification}, \orgname{Sandia National Laboratories}, \orgaddress{ \city{Albuquerque}, \postcode{87185-1320}, \state{New Mexico}, \country{USA}}}

\affil[2]{\small \orgdiv{Scientific Machine Learning}, \orgname{Sandia National Laboratories},\\ \orgaddress{ \city{Albuquerque}, \postcode{87185-1320}, \state{New Mexico}, \country{USA}}}



\abstract{
We develop a parallel-in-time multigrid preconditioner for augmented systems.
These saddle-point systems are foundational to numerical optimization.
Our preconditioner, when paired with a suitable optimization method,
accelerates the solution of optimal control problems.
We construct the preconditioner by introducing virtual interface variables
that enable time-domain decomposition.
After permuting the resulting augmented system into block tridiagonal form,
we develop a geometric multigrid scheme with a block Jacobi smoother, which
parallelizes trivially in time.
As the coarse grid solver we use GMRES preconditioned with a symmetric
Gauss-Seidel iteration.
We use the multigrid scheme to precondition a flexible GMRES~\cite{saad.1993}
iteration for the solution of the augmented system.
We combine our preconditioner with the matrix-free sequential quadratic
programming~(SQP) algorithm \cite{heinkenschloss.2014} to solve optimal
control problems involving the van der Pol oscillator and the
viscous Burgers' equation.
We find that the preconditioner is remarkably effective when the problems
are suitably scaled.
}



\maketitle

%
%
%
\section{Introduction}
\label{sec:intro}

This paper is motivated by optimal control of differential equations.
In particular, we study the optimization problem
\begin{subequations}
\label{eq:optc-aj}
\begin{align}
  \underset{\mathcal{u}\in\mathcal{U},\;\mathcal{z}\in\mathcal{Z}}{\text{minimize}}&\quad
\mathcal{J}(\mathcal{u},\mathcal{z}) \\
\text{subject to}&\quad \mathcal{c}(\mathcal{u},\mathcal{z}) := M\mathcal{u} - \int_0^t f(\mathcal{u},\mathcal{z}) = 0.
\label{eq:ode-aj}
\end{align}
\end{subequations}
Instances of \eqref{eq:optc-aj} can be used to enforce additional constraints,
be they equalities or inequalities, using, e.g., augmented Lagrangian
methods \cite{nocedal.2006,antil.2023}.
The domains~$\mathcal{U}$ and~$\mathcal{Z}$ correspond to a state space and
a control space, respectively.
For concreteness, let
$$\mathcal{U}=L^2([0,1];\mathbb{R}^p) \quad\text{and}\quad \mathcal{Z}=L^2([0,1];\mathbb{R}^q).$$
The objective function $\mathcal{J}:\mathcal{U}\times\mathcal{Z}\to\mathbb{R}$ is a
smooth measure of cost, and---as is conventional \cite{malek.2014}---the constraint
function $\mathcal{c}:\mathcal{U}\times\mathcal{Z}\to\mathcal{U}^*$, where
$\mathcal{U}^*$ is the topological dual space of $\mathcal{U}$.
Within the constraint function, $M:\mathcal{U}\to\mathcal{U}^*$ is linear
and $f:\mathcal{U}\times\mathcal{Z}\to\mathcal{U^*}$.
Note that \eqref{eq:ode-aj} is a differential equation
for~$\mathcal{u}$ given~$\mathcal{z}$.
For sufficiently regular $\mathcal{u}$, \eqref{eq:ode-aj} is equivalent to
an initial value problem on the time domain $[0,1]$, i.e.,
\begin{align*}
M \frac{d\mathcal{u}}{dt} = f(\mathcal{u},\mathcal{z}), \quad u(0) = 0.\footnotemark
\end{align*}
\footnotetext{
This homogeneous initial condition is without loss of generality since a
nonzero initial value $u_0$ can be absorbed into $\mathcal{u}$ with the
shift $\mathcal{u}\gets\mathcal{u}+\chi$, where $\chi$ is the function
on $[0,1]$ everywhere equal to $u_0$.
Similarly, an invertible $M$ can be absorbed into $f$ with $f\gets Mf$.
We avoid this reduction, however, since it can ($i$) destroy sparsity
inherent to the problem and ($ii$) change the co-domain of the residual,
complicating the presentation of our work.
}

The state of practice for solving \eqref{eq:optc-aj} at scale is to use
derivative-based numerical optimization.
When doing so, it is common to reformulate the problem \eqref{eq:optc-aj}
posed over $\mathcal{U}\times\mathcal{Z}$ into an equivalent problem
over $\mathcal{Z}$ only.
Concretely, if \eqref{eq:ode-aj} admits a solution
operator~$\mathcal{S}:\mathcal{Z}\to\mathcal{U}$, then \eqref{eq:ode-aj}
can be enforced implicitly, yielding the reduced problem
\begin{align}\label{eq:red-aj}
  & \underset{\mathcal{z}\in\mathcal{Z}}{\text{minimize}}\;\; \mathcal{J}(\mathcal{S}(\mathcal{z}),\mathcal{z}).
\end{align}

\begin{figure}
\centering
\includegraphics[width=0.9\textwidth]{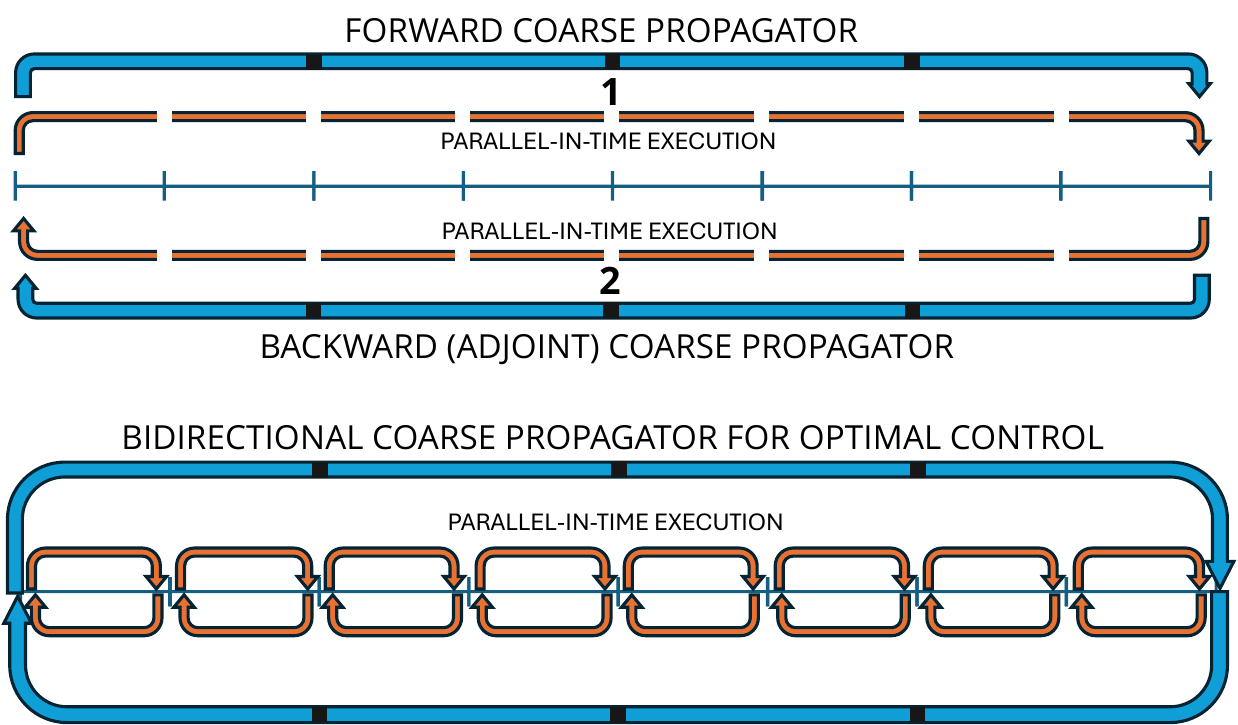}
\caption{
Schematic of parallel-in-time integration applied to optimal control.
Black squares indicate serialization points, belonging to the coarse
correction/propagation schemes.
\emph{Top:}
Conventional approaches remain intrinsically serial, consisting of forward
time stepping, 1, to evaluate the objective function, followed by backward
time stepping, 2, to evaluate derivatives.
\emph{Bottom:}
Our preconditioner is, in contrast, a parallel-in-time method tailored to
the coupled bidirectional flow of time in optimality systems arising
from optimal control.
Parallel-in-time execution involves concurrent solutions of optimal control
problems on the time subdomains.}
\label{fig:parallel}
\end{figure}

Let $S$ be a discretization of $\mathcal{S}$.
We address the case where the evaluation of $S$ is a matter of forward time
stepping and the efficient evaluation of its derivative is a matter of
backward time stepping to solve an adjoint equation.
The serial nature of this time stepping can be a computational bottleneck
that results in long solution times, even for simple problems.
The parallel-in-time methods~\cite{maday.2002,ulbrich.2007,guenther.2019,guenther.2020}
ameliorate this time stepping bottleneck but retain a key serial attribute:
they apply parallel-in-time integration \cite{lions.2001,gander.2015} twice,
in sequence, first to evaluate~$S$ and then to evaluate its derivative.
See Figure~\ref{fig:parallel}.

Our contribution is a parallel-in-time multigrid preconditioner that avoids
this serial structure for a specific yet foundational instance of optimal
control.
The preconditioner pertains to \emph{augmented systems},
\begin{align}\label{eq:augsys-aj}
&Ax = b, \quad A = \begin{bmatrix} I & B^\top \\ B & 0 \end{bmatrix},
\end{align}
which we take to be finite-dimensional.
In particular, $B$ is a matrix with at least as many columns as rows
and~$I$ is an identity of the appropriate dimension.
The quantities~$x$ and~$b$ are vectors with dimension equal to the number
of columns plus the number of rows of $B$.
Augmented systems are saddle-point systems~\cite{benzi.2005}.
They include the necessary optimality conditions for discretizations of the
linear-quadratic instance of \eqref{eq:optc-aj},
\begin{subequations}
\label{eq:tracking-aj}
\begin{align}
\underset{\mathcal{u}\in\mathcal{U},\;\mathcal{z}\in\mathcal{Z}}{\text{minimize}}&\quad
  \frac{1}{2} \| \mathcal{u} - \mathcal{\hat{u}} \|^2
+ \frac{1}{2} \| \mathcal{z} - \mathcal{\hat{z}} \|^2
\label{eq:tracking1-aj}\\
\text{subject to}&\quad M\mathcal{u} - \int_0^t (B_1\mathcal{u} + B_2\mathcal{z} + g) = 0,
\label{eq:tracking2-aj}
\end{align}
\end{subequations}
where $\mathcal{\hat{u}}\in\mathcal{U}$, $\mathcal{\hat{z}}\in\mathcal{Z}$,
$B_1:\mathcal{U}\to\mathcal{U}^*$ and $B_2:\mathcal{Z}\to\mathcal{U}^*$ are
linear, and $g\in\mathcal{U}^*$.
Augmented systems are also, e.g., the key computational kernel of
implementations of the sequential quadratic programming (SQP) optimization
algorithm \cite{heinkenschloss.2014} that can be used to solve the general
problem \eqref{eq:optc-aj}.

Other bidirectional parallel-in-time methods for optimal control have
been proposed; see \cite{gander.2020,lin.2022,cyr.2024}.
To our knowledge, our preconditioner is ($i$) the first that preserves the
full saddle-point structure of optimality systems---enabling block scalings
that considerably accelerate convergence of the Krylov scheme---and
($ii$) most closely connected with multigrid, supporting standard restriction
and interpolation operators as well as cycles of arbitrary type and depth.

The remainder of the paper is organized as follows.
In Section~\ref{sec:notation-aj}, we introduce notation.
In Section~\ref{sec:augsys-aj}, we discuss the significance of augmented
systems.
We show that they include least-squares and minimum-norm problems as special
cases, and we explain their connection with \eqref{eq:tracking-aj}.
In Section~\ref{sec:disc-aj}, we discretize~\eqref{eq:ode-aj} to obtain
augmented systems that converge to their infinite-dimensional analogs
as the length of the time step approaches zero.
We develop our parallel-in-time multigrid preconditioner in
Section~\ref{sec:Aug_pint}, including the appropriate time-domain
decomposition, smoother, interpolation and restriction, and coarse-grid solver. 
In Section~\ref{sec:numerics}, we present numerical results.
We combine the parallel-in-time multigrid preconditioner
with~\cite{heinkenschloss.2014} to solve \eqref{eq:optc-aj}
when~\eqref{eq:ode-aj} is the van der Pol oscillator and the viscous
Burgers' equation.
Finally, we demonstrate that our preconditioner can also be used to
efficiently solve the inviscid Burgers' equation---a hyperbolic nonlinear
equation.

%
%
%
\section{Notation}
\label{sec:notation-aj}

We use the symbol $\partial$ for the Jacobian of a function with
respect to its inputs and a dot for the special case of the total
time derivative, e.g., $\dot{\mathcal{x}} := \frac{d\mathcal{x}}{dt}$.
We add a subscript to $\partial$ when the Jacobian is with respect
to a subset of inputs.
As is standard, $\mathbb{R}^d$ is the set of real-valued $d$-tuples,
though we sometimes overload notation slightly, when the context is clear,
by taking $\mathbb{R}^d$ to also mean $d$-dimensional Euclidean space
(the Hilbert space of these tuples with the usual inner product).

%
%
%
\section{Augmented Systems}
\label{sec:augsys-aj}

The foundational importance of augmented systems in optimization
arises from their being the necessary optimality conditions for
least-squares problems involving $B^\top$ and minimum norm problems
involving $B$.
These necessary conditions are straightforward to derive.
The least-squares problem
$$\underset{u}{\text{minimize}}\;\; \frac{1}{2}\|\chi - B^\top u\|^2 ,$$
for suitable vectors $u$ and $\chi$, is dual to the orthogonal projection
on the null space of $B$, i.e.,
\begin{align}\label{eq:ls-aj}
\begin{split}
  \underset{s}{\text{minimize}}&\quad\frac{1}{2}\|\chi - s\|^2 \\
\text{subject to}&\quad Bs = 0.
\end{split}
\end{align}
The orthogonal projection problem \eqref{eq:ls-aj} has the Lagrangian
\begin{align*}
&\frac{1}{2}\|\chi - s\|^2 + \lambda^\top Bs,
\end{align*}
and requiring its gradients with respect to $s$ and $\lambda$ to be zero
results in an instance of \eqref{eq:augsys-aj} with
\begin{align}\label{eq:augsysls-aj}
x = \begin{bmatrix}
s \\ \lambda
\end{bmatrix}
\quad\text{and}\quad
b = \begin{bmatrix}
\chi \\ 0
\end{bmatrix}.
\end{align}
Meanwhile, the minimum norm problem is
\begin{align}\label{eq:mn-aj}
\begin{split}
  \underset{s}{\text{minimize}}&\quad \frac{1}{2}\|s\|^2 \\
  \text{subject to}&\quad Bs = \chi.
\end{split}
\end{align}
In this setting, $\chi$ appears in the gradient with respect to $\lambda$,
as opposed to $s$.
Thus, an instance of \eqref{eq:augsys-aj} again results as the necessary
optimality conditions, with $x$ as in \eqref{eq:augsysls-aj} and
$$b = \begin{bmatrix}
0 \\ \chi
\end{bmatrix}.$$

In the next section, we discretize \eqref{eq:ode-aj} into an algebraic
equation $F(u,z) = 0$ such that augmented systems with $B = \partial F$ are
consistent with their infinite-dimensional analogs.
These analogs are linear equations posed over function spaces and include the
necessary optimality conditions for \eqref{eq:tracking-aj}.
The necessary optimality conditions for that problem are readily derived in a
manner similar to those for \eqref{eq:ls-aj} and \eqref{eq:mn-aj}.
See, e.g., \cite{clarke.2013}.
They are a boundary value problem (BVP) involving the Lagrange
multiplier~$\Lambda\in\mathcal{U}^{**}=\mathcal{U}$ for \eqref{eq:tracking2-aj},
\begin{align}\label{eq:bvp1-aj}
\begin{split}
\mathcal{u} + M^\top\Lambda - B_1^\top\int_t^1 \Lambda &= \hat{\mathcal{u}}  \\
\mathcal{z} - B_2^\top\int_t^1 \Lambda &= \hat{\mathcal{z}} \\
\quad M\mathcal{u} - \int_0^t (B_1\mathcal{u} + B_2\mathcal{z}) &= \int_0^t g.
\end{split}
\end{align}
Note that this BVP has a saddle-point structure similar to \eqref{eq:augsys-aj}.
We will see from $F$ that---as one might intuitively expect---\eqref{eq:ls-aj}
corresponds to the forcing term $g$ being zero and \eqref{eq:mn-aj} corresponds
to the tracking terms $\hat{\mathcal{u}}$ and $\hat{\mathcal{z}}$ being zero.

\subsection{A Case for Multigrid in Time}

Let $\Xi(t) := \int_t^1\Lambda$ and suppose $\Lambda$ and $\mathcal{u}$ are regular enough for the second fundamental theorem of calculus to hold, in which case we can write \eqref{eq:bvp1-aj} as
\begin{align}\label{eq:bvp2-aj}
\begin{split}
  \mathcal{u} - M^\top\dot{\Xi} - B_1^\top\Xi &= \hat{\mathcal{u}}  \\
  \mathcal{z} - B_2^\top\Xi &= \hat{\mathcal{z}}  \\
  M\dot{\mathcal{u}} - (B_1\mathcal{u} + B_2\mathcal{z}) &= g.
\end{split}
\end{align}
Using the first and second equations of \eqref{eq:bvp2-aj} to eliminate $\mathcal{u}$ and $\mathcal{z}$ yields the second-order differential equation
\begin{align}
\label{eq:ellipticPDE}
-\ddot{\Xi} + (B_1-B_1^T)\dot{\Xi} + (B_1 B_1^\top + B_2 B_2^\top) \Xi
= g + \dot{\hat{\mathcal{u}}} - B_1\hat{\mathcal{u}} - B_2\hat{\mathcal{z}}.
\end{align}
This differential equation is the Schur complement of the $\mathcal{u}$ and $\mathcal{z}$ identity block in \eqref{eq:bvp2-aj}.
Importantly, \eqref{eq:ellipticPDE} resembles an \emph{elliptic} equation in
time, suggesting multigrid (in time) for augmented systems is more natural
than multigrid for linearizations of the forward equation \eqref{eq:ode-aj}.
This observation is not without precedent.
The necessary optimality conditions for tracking problems with parabolic
partial differential equations (as opposed to the ordinary differential
equations we consider) are elliptic \cite{cyr.2024,gander.2016}, and reduced
Hessians for tracking problems have been shown to sometimes exhibit elliptic
behavior even when the differential equation is hyperbolic \cite{lewis.2005}.

%
%
%
\section{Discretization}
\label{sec:disc-aj}

We discretize \eqref{eq:ode-aj} with a time-stepping scheme, which we take
to be the $\theta$-method with a piecewise constant control; namely,
\begin{multline}\label{eq:theta-aj}
  c_k := \bM \uu_k - \bM \uu_{k-1} - \Delta t \left[\theta f(\uu_k,\zz_k)  + (1 - \theta) f(\uu_{k-1},\zz_k)\right] = 0, \\
   k = 1, 2, \ldots, n, \quad \theta\in (0,1).
\end{multline}
Here, $u_k\in\mathbb{R}^p$ approximates the state~$\mathcal{u}$ at~$k\Delta t$
for~$\Delta t = 1/n$ and~$u_0 := 0$.
Analogously, $z_k\in\mathbb{R}^q$ approximates the control~$\mathcal{z}$ over
the interval~$$[(k-1)\Delta t,\; k\Delta t].$$
The discretized state and control variables are thus
$$u = \begin{bmatrix} u_1 \\ u_2 \\ \vdots \\ u_n \end{bmatrix}\in\mathbb{R}^{pn} \quad\text{and}\quad
  z = \begin{bmatrix} z_1 \\ z_2 \\ \vdots \\ z_n \end{bmatrix}\in\mathbb{R}^{qn}. $$
We think of the discretized constraint function as vector valued:
$$\begin{bmatrix} c_1 \\ c_2 \\ \vdots \\ c_n \end{bmatrix} =: c(u,z) \in\mathbb{R}^{pn}.$$

We note that the time-stepping scheme \eqref{eq:theta-aj} is in the interest of concreteness.
Our parallel-in-time multigrid preconditioner supports other discretizations, to include alternative approximations of the state like fourth order Runge-Kutta and alternative approximations of the control.

\subsection{Scaling}\label{sec:scaling-aj}

Care is required to ensure augmented systems obtained from the linearization
of~\eqref{eq:theta-aj} are consistent with their infinite dimensional
analog \eqref{eq:bvp1-aj}.
The domain of the continuous constraint function~$\mathcal{c}$ is the product
space~$\mathcal{U}\times\mathcal{Z}$.
Recall that $\mathcal{U}=L^2([0,1];\mathbb{R}^p)$, which has the inner product
$$(\mathcal{u},\tilde{\mathcal{u}})_{\mathcal{U}} := \int_0^1 \tilde{\mathcal{u}}^\top\mathcal{u}
\quad\text{for all}\quad \mathcal{u},\tilde{\mathcal{u}}\in\mathcal{U}.$$
A natural choice for the discrete state space $U$ is
therefore~$\mathbb{R}^{pn}$ equipped with the inner product
\begin{align}\label{eq:ipu-aj}
&(u,\tilde{u})_U = \tilde{u}^\top \Sigma_U u
\quad\text{for all}\quad u, \tilde{u}\in U,
\end{align}
where $\Sigma_U$ is the $pn\times pn$ identity scaled by $\Delta t$.
Analogously, we choose the discrete control space $Z$ to
be~$\mathbb{R}^{qn}$ equipped with the inner product
\begin{align}\label{eq:ipz-aj}
&(z,\tilde{z})_Z = \tilde{z}^\top \Sigma_Z z
\quad\text{for all}\quad z, \tilde{z}\in Z,
\end{align}
where $\Sigma_Z$ is the $qn\times qn$ identity scaled by $\Delta t$.
In a slight abuse of notation, we let $u$ and $z$ be the state and control
components of a vector~$x\in U\times Z$.
As the inner product on~$U\times Z$, we choose
\begin{align}\label{eq:ipx-aj}
&(x,\tilde{x})_X = (u,\tilde{u})_U + \gamma^2(z,\tilde{z})_Z
\quad\text{for all}\quad x, \tilde{x}\in U\times Z,
\end{align}
where $\gamma$ is a positive scalar that will be key to
designing a performant multigrid scheme. 
Recall that the codomain of $\mathcal{c}$ is $\mathcal{U}^*$.
We state the following straightforward fact about the discretization of this space.
\begin{proposition}\label{p:pd}
  Let $U^*$ be $\mathbb{R}^{pn}$ equipped with the inner product
\begin{align}\label{eq:ipud-aj}
&(u^*,\tilde{u}^*)_{U^*} = (\tilde{u}^*)^\top \Sigma_{U^*} u^*
\quad\text{for all}\quad u^*, \tilde{u}^*\in U^*
\end{align}
for some matrix $\Sigma_{U^*}$.
If the primal-dual pairing $\langle\cdot,\cdot\rangle_{U,U^*}$ is the Euclidean inner product between vectors in $\mathbb{R}^{pn}$, then $\Sigma_{U^*} := \Sigma_U^{-1}$
\begin{itemize}\small
\item
induces the operator norm on $U^*$
\item
is an isomorphism taking $U^*$ into $U$.
\end{itemize}
\end{proposition}
We define $U^*$ in accordance with Proposition~\ref{p:pd}.
Just as $\mathcal{c}:\mathcal{U}\times\mathcal{Z}\to\mathcal{U}^*$, we have $$c:U\times Z\to U^*.$$
As a result, $\partial c$ evaluated at a point in $U\times Z$ is a linear map from $U\times Z$ to $U^*$, meaning augmented systems involving $\partial c$ would have an output space dual to their input space, necessitating a non-identity $(1,1)$ block in general.
Complications like this are avoided by ``transplanting'' the domain and codomain of $c$ to Euclidean space.
A matrix $\Sigma$ inducing an inner product must be symmetric positive definite and hence has a square root that we denote as $\Sigma^{\frac{1}{2}}$ \cite{golub.2013}.
This square root defines an isomorphism from the associated inner product space into Euclidean space, which is how we transplant.
In particular, we take as our discretization of $\mathcal{c}$ not $c$ but the operator
$$
F:\mathbb{R}^{pn}\times\mathbb{R}^{qn}\to\mathbb{R}^{pn}
\quad\text{with}\quad
F(u,z):=\Sigma_{U^*}^{\;\;\frac{1}{2}}c\left(\,\Sigma_U^{-\frac{1}{2}}u,\,(\gamma^2\Sigma_Z)^{-\frac{1}{2}}z\,\right).
$$

Having arrived at a discretization of \eqref{eq:ode-aj}, we present in the next section our parallel-in-time multigrid preconditioner for augmented systems with $B=\partial F$.
To do so, we find it helpful to label the different pieces of $\partial F$.
For $k=1,2,\ldots, n$, let
\begin{align}\label{eq:blocks-aj}
\begin{split}
N_k
&= \phantom{(\gamma)}\Delta t^{-1} \Big[\hspace{1.2em} \bM - \hspace{2.5em}\theta\Delta t \partial_u f\left(\,\Delta t^{-\frac{1}{2}}\uu_{k\phantom{-1}},\,(\gamma^2\Delta t)^{-\frac{1}{2}}\zz_k\,\right) \Big]  \\
O_k
&= \phantom{(\gamma)}\Delta t^{-1}\Big[- \bM - (1-\theta)\Delta t \partial_u f\left(\,\Delta t^{-\frac{1}{2}}\uu_{k-1},\,(\gamma^2\Delta t)^{-\frac{1}{2}}\zz_k\,\right) \Big]  \\
Z_k
&= (\gamma\Delta t)^{-1} \Big[
\hspace{1.2em}\phantom{M}-\hspace{2.6em}\theta \Delta t \partial_z f\left(\,\Delta t^{-\frac{1}{2}}\uu_{k\phantom{-1}},\,(\gamma^2\Delta t)^{-\frac{1}{2}}\zz_k\,\right)
\\
&\hspace{7.7em}- (1-\theta) \Delta t\partial_z f\left(\,\Delta t^{-\frac{1}{2}}\uu_{k-1},\,(\gamma^2\Delta t)^{-\frac{1}{2}}\zz_k\,\right) \Big].
\end{split}
\end{align}

%
%
%
\section{Preconditioner} \label{sec:Aug_pint}

Figure~\ref{matrix:eq:long:serial} shows the algebraic structure of the augmented system \eqref{eq:augsys-aj} when $B=\partial F$.
In keeping with Section~\ref{sec:augsys-aj}, we think of $\lambda$ as the adjoint variable associated with $B$.
As such, $\lambda\in\mathbb{R}^{pn}$ and we take $\lambda_k$ to be the $p$ components of $\lambda$ associated with \eqref{eq:theta-aj} at the subscripted value of $k$.
\begin{figure}[h]
\begin{align*}
\scalemath{0.7}{
\begin{pNiceArray}{c c c c |  c c c c  | c c c c }
\textcolor{green}{I} &  &  &  &  &  &  &  & {\textcolor{red}{ \bK_1^\TT}} & \textcolor{orange}{\bC_2^\TT} &  &     \\
 & \textcolor{green}{I} &  &  &  &  &  &  &  & \textcolor{red}{ \bK_2^\TT} & \textcolor{orange}{\bC_3^\TT} &  \\
 &  & \textcolor{green}{I}  &  &  &  &  &  &  &  & \textcolor{red}{ \bK_3^\TT} &  \textcolor{orange}{\bC_4^\TT} \\
 &  &   & \textcolor{green}{I} &  &  &  &  &  &  & & \textcolor{red}{ \bK_4^\TT} \\
\arrayrulecolor[rgb]{0.8,0.8,0.8} \hline
 &  &   &  &  \textcolor{green}{I} &  &  &  & \textcolor{blue}{\bB_1^\TT} &  & &  \\
 &  &   &  &  &  \textcolor{green}{I} &  &  &  & \textcolor{blue}{\bB_2^\TT} & &  \\
 &  &   &  &  &  &  \textcolor{green}{I} &  &  &  &\textcolor{blue}{\bB_3^\TT}  &  \\
 &  &   &  &  &  &  &  \textcolor{green}{I} &  &  & & \textcolor{blue}{\bB_4^\TT} \\
\hline
\textcolor{red}{ \bK_1} & &  &  & \textcolor{blue}{\bB_1} &  &  &  &  &  & &  \\
 \textcolor{orange}{\bC_2} & \textcolor{red}{\bK_2} &   &  &  &\textcolor{blue}{\bB_2}   &  &  &  &  & &  \\
 &  \textcolor{orange}{\bC_3} & \textcolor{red}{\bK_3}  &  &  &  & \textcolor{blue}{\bB_3} &  &  &  & &  \\
 &  &  \textcolor{orange}{\bC_4}  & \textcolor{red}{\bK_4} &  &  &  & \textcolor{blue}{\bB_4}
\end{pNiceArray}
\begin{pNiceArray}{c}
\uu_1\\
\uu_2\\
\uu_3\\
\uu_4\\
\arrayrulecolor[rgb]{0.8,0.8,0.8}
\hline
\zz_1\\
\zz_2\\
\zz_3\\
\zz_4\\
\hline
\bla_1\\
\bla_2\\
\bla_3\\
\bla_4
\end{pNiceArray}
=
\begin{pNiceArray}{c}
\\[1.5em]
\uparrow\\
\\[1.5em]
b\\
\\[1.5em]
\downarrow\\
\\[1.5em]
\end{pNiceArray}
}
\end{align*}
\caption{Augmented system when $n = 4$.}
\label{matrix:eq:long:serial}
\end{figure}

Analogous to \eqref{eq:bvp1-aj}, the lower triangular structure
of~$\partial_u F$, block (3,1) in~Figure~\ref{matrix:eq:long:serial},
flows information forward in time, while its transpose, block (1,3)
in~Figure~\ref{matrix:eq:long:serial}), flows information backward in time.
We exploit this bidirectionality by adding virtual variables at the
interfaces of time steps, as is typical in domain decomposition
methods for spatial discretizations \cite{smith.2004}.

\subsection{Virtual Variables}
\label{sec:lift}

\begin{figure}[!htb]
\centering
\includegraphics[width=.675\textwidth]{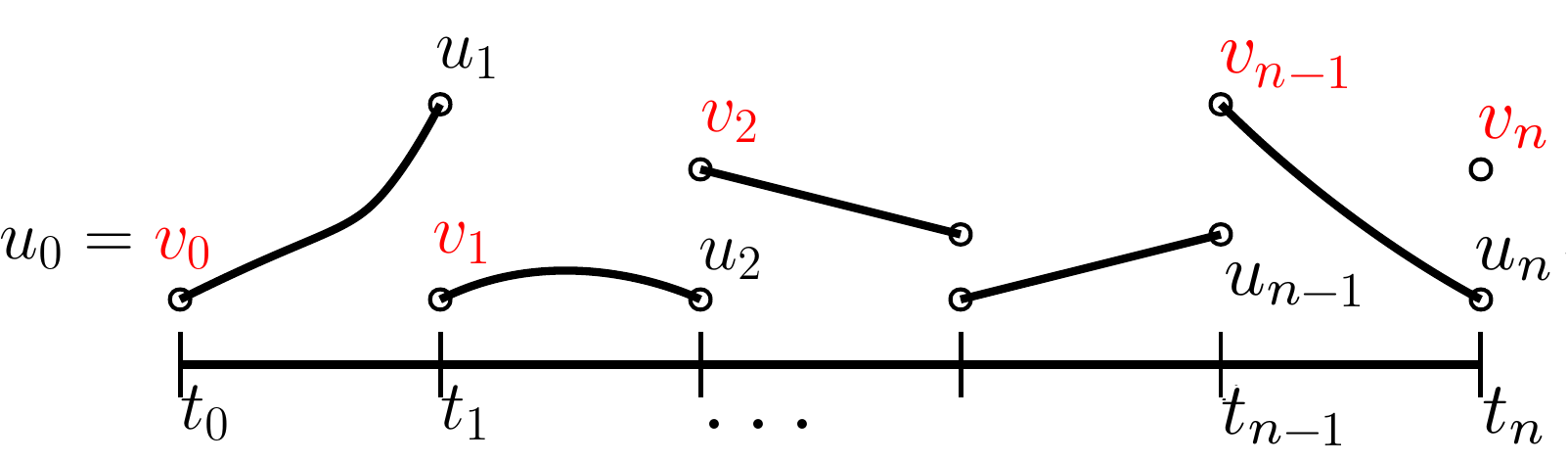}
\caption{Time domain decomposition.}
\label{fig:relax}
\end{figure}

We introduce virtual interface variables $v\in\mathbb{R}^{pn}$ so that in place of \eqref{eq:theta-aj} we have
\begin{multline}\label{eq:theta2-aj}
  \bM \uu_k - \bM \vv_{k-1} - \Delta t \left[\theta f(\uu_k,\zz_k)  + (1 - \theta) f(\vv_{k-1},\zz_k)\right] = 0, \\
   k = 1, 2, \ldots, n, \quad \theta\in (0,1).
\end{multline}
Here, we define $v_k$ in a similar manner as $u_k$ with $v_0 := 0$.
To the constraints \eqref{eq:theta2-aj}, we add
\begin{align}\label{eq:equal-aj}
  v = u,
\end{align}
to obtain an approximation equivalent to \eqref{eq:theta-aj}.
Figure~\ref{fig:relax} shows the effect of this change.
When $u\neq v$, as is generally the case except at the solution,
the virtual variables promote decoupling of the time steps,
i.e., enable concurrent local computations.

Analogous to \eqref{eq:theta-aj}, we interpret \eqref{eq:theta2-aj} as prescribing a residual in $U^*$.
We formulate \eqref{eq:equal-aj} in the same way by writing the constraint as
$$d(u,v) := \Delta t(u-v) = 0.$$
The scaling considerations in Section~\ref{sec:scaling-aj} imply that $d$ transplants into the Euclidean space function
$$
\Sigma_{U^*}^{\;\;\frac{1}{2}}d\left(\,\Sigma_U^{-\frac{1}{2}}u,\,\Sigma_U^{-\frac{1}{2}}v\,\right) = u-v.
$$
Let the Euclidean space function
$$\tilde{F}:\mathbb{R}^{pn}\times\mathbb{R}^{pn}\times\mathbb{R}^{qn}
\to\mathbb{R}^{pn}\times\mathbb{R}^{pn}$$
be the scaled discretization for \eqref{eq:theta2-aj} and \eqref{eq:equal-aj}, i.e., the counterpart of $F$ when introducing virtual variables.

Figure~\ref{matrix:eq:long:v} shows the algebraic structure of the augmented system \eqref{eq:augsys-aj} when $B=\partial \tilde{F}$.
The blocks $V_k$ and $\tilde{Z}_k$ are $O_k$ and $Z_k$, respectively, but with $u_{k-1}$ replaced by $v_{k-1}$; that is,
\begin{align*}
V_k
&= \phantom{(\gamma)}\Delta t^{-1}\Big[- \bM - (1-\theta)\Delta t \partial_u f\left(\,\Delta t^{-\frac{1}{2}}\vv_{k-1},\,(\gamma^2\Delta t)^{-\frac{1}{2}}\zz_k\,\right) \Big]
\\
\tilde{Z}_k
&= (\gamma\Delta t)^{-1} \Big[
\hspace{1.2em}\phantom{M}-\hspace{2.6em}\theta \Delta t \partial_z f\left(\,\Delta t^{-\frac{1}{2}}u_{k\phantom{-1}},\,(\gamma^2\Delta t)^{-\frac{1}{2}}\zz_k\,\right)
\\
&\hspace{7.7em}- (1-\theta) \Delta t\partial_z f\left(\,\Delta t^{-\frac{1}{2}}v_{k-1},\,(\gamma^2\Delta t)^{-\frac{1}{2}}\zz_k\,\right) \Big].
\end{align*}
The adjoint variable $\lambda$ is now associated with \eqref{eq:theta2-aj} instead of \eqref{eq:theta-aj}, and $\mu\in\mathbb{R}^{pn}$ is the adjoint variable associated with \eqref{eq:equal-aj}.
We define $\mu_k$ in a manner similar to $\lambda_k$---it is the $p$ components of $\mu$ corresponding to the equality of $u_k$ and $v_k$.
\begin{figure}[h]
\begin{align*}
\scalemath{0.65}{
\begin{pNiceArray}{cccc|cccc|cccc|cccc|cccc}
\textcolor{green}{I} & & & & & & & & & & & & \textcolor{red}{N_1^\top} & & & & \phantom{-}I
\\
& \textcolor{green}{I} & & & & & & & & & & & & \textcolor{red}{N_2^\top} & & & & \phantom{-}I
\\
& & \textcolor{green}{I} & & & & & & & & & & & & \textcolor{red}{N_3^\top} & & & & \phantom{-}I
\\
& & & \textcolor{green}{I} & & & & & & & & & & & & \textcolor{red}{N_4^\top} & & & & \phantom{-}I
\\
\arrayrulecolor[rgb]{0.8,0.8,0.8}
\hline
& & & & \textcolor{green}{I} & & & & & & & & & \textcolor{orange}{V_2^\top} & & & -I
\\
& & & & & \textcolor{green}{I} & & & & & & & & & \textcolor{orange}{V_3^\top} & & & -I
\\
& & & & & & \textcolor{green}{I} & & & & & & & & & \textcolor{orange}{V_4^\top} & & & -I
\\
& & & & & & & \textcolor{green}{I} & & & & & & & & & \phantom{V_4^\top} & & & -I
\\
\hline
& & & & & & & & \phantom{-}\textcolor{green}{I} & & & & \textcolor{blue}{\tilde{Z}_1^\top}
\\
& & & & & & & & & \phantom{-}\textcolor{green}{I} & & & & \textcolor{blue}{\tilde{Z}_2^\top}
\\
& & & & & & & & & & \phantom{-}\textcolor{green}{I} & & & & \textcolor{blue}{\tilde{Z}_3^\top}
\\
& & & & & & & & & & & \phantom{-}\textcolor{green}{I} & & & & \textcolor{blue}{\tilde{Z}_4^\top}
& \phantom{N_1} & \phantom{N_2} & \phantom{N_3} & \phantom{N_4}
\\
\hline
 \textcolor{red}{N_1} &     &     &     \phantom{-}\textcolor{orange}{\phantom{V_1}} &     &     &     & & \textcolor{blue}{\tilde{Z}_1} &       &       &       \\
     & \textcolor{red}{N_2} &     &                   & \phantom{-}\textcolor{orange}{V_2} &     &     & &       & \textcolor{blue}{\tilde{Z}_2} &       &       \\
     &     & \textcolor{red}{N_3} &                   &     & \phantom{-}\textcolor{orange}{V_3} &     & &       &       & \textcolor{blue}{\tilde{Z}_3} &       \\
     &     &     & \textcolor{red}{N_4}               &     &     & \phantom{-}\textcolor{orange}{V_4} & &       &       &       & \textcolor{blue}{\tilde{Z}_4} \\
\hline
\phantom{-}I & & & & -I       \\
& \phantom{-}I & & & & -I     \\
& & \phantom{-}I & & & & -I   \\
& & & \phantom{-}I & & & & -I
\end{pNiceArray}
\begin{pNiceArray}{c}
\uu_1\\
\uu_2\\
\uu_3\\
\uu_4\\
\arrayrulecolor[rgb]{0.8,0.8,0.8}
\hline
\vv_1\\
\vv_2\\
\vv_3\\
\vv_4\\
\hline
\zz_1\\
\zz_2\\
\zz_3\\
\zz_4\\
\hline
\bla_1\\
\bla_2\\
\bla_3\\
\bla_4\\
\hline
\mmu_1\\
\mmu_2\\
\mmu_3\\
\mmu_4
\end{pNiceArray}
=
\begin{pNiceArray}{c}
\\[4em]
\uparrow\\
\\[4em]
b\\
\\[4em]
\downarrow\\
\\[4em]
\end{pNiceArray}
}
\end{align*}
\caption{Augmented system with virtual variables when $n = 4$.}
\label{matrix:eq:long:v}
\end{figure}

\subsection{Permutation of the Augmented System}
\label{sub:section:reordering}

The augmented system in~Figure~\ref{matrix:eq:long:v} is organized to
have the variable and constraint types changing more slowly than the
time indices $k$.
For instance, the vector to the immediate right of the matrix
in~Figure~\ref{matrix:eq:long:v} consists of the states at all times,
followed by the interface variables at all times, and so forth.
We call this \emph{type-major} order.
To develop our parallel-in-time multigrid preconditioner,
we permute the augmented system into \emph{time-major} order,
as shown in Figure~\ref{matrix:eq:ord:long}.
This change is a matter of replacing $Ax=b$ in \eqref{eq:augsys-aj}
with the equivalent system
\begin{align}\label{eq:augsysp-aj}
  &(PAQ^\top)(Qx)=Pb
\end{align}
for suitable permutation matrices $P$ and $Q$.
We omit algebraic expressions for these permutation matrices since
they are unwieldy and not needed in practice.
They are straightforward to obtain, however, using shuffle matrices
and Kronecker products~\cite{golub.2013}.
We define $\tilde{A} = PAQ^\top$, $y = Qx$,
and $\tilde{b} = Pb$ so that \eqref{eq:augsysp-aj} becomes
\begin{align}\label{eq:augsysp2-aj}
  &\tilde{A}y = \tilde{b}.
\end{align}
\begin{figure}[h]
\begin{align*}
\scalemath{0.7}{
\begin{pNiceArray}{ccc|cccccc|ccccc|ccccc|cc}
\textcolor{green}{I}  &  &\textcolor{red}{\bK_1^\top}&  & & & & & \phantom{-}I & & & & & & & & & & & &\\
 & \textcolor{green}{I}& \textcolor{blue}{\bB_1^\top} & & & & & & &  & & & & & & & & & & &\\
 \textcolor{red}{\bK_1}  & \textcolor{blue}{\bB_1} & & & & & & & & & & & & & & & & & & &\\
\arrayrulecolor[rgb]{0.8,0.8,0.8} \hline
 & & & & \textcolor{green}{I}& & & \textcolor{orange}{V_1^\top} & -I  & & & & & & & & & & & &\\
 & & & & &\textcolor{green}{I} & & \textcolor{red}{\bK_2^\top}& & & & &  &\phantom{-}I  &  & & & & & &\\
 & & & & & & \textcolor{green}{I} & \textcolor{blue}{\tilde{Z}_2^\top} & & & & & & & & & & & & &\\
 & & & & \textcolor{orange}{V_2} & \textcolor{red}{N_2} &\textcolor{blue}{\tilde{Z}_2} & & & & & & & & & & & & & &\\
I& & &  & -I & & & &  & & & & & & & & & & & &\\
 \hline
 &  & & & & & & & & \phantom{-}\textcolor{green}{I} & & & \textcolor{orange}{V_3^\TT}& -I  & & & & & & &\\
 & & &  & & & & & & & \textcolor{green}{I}  & &  \textcolor{red}{N_3^\top} & & & & & & \phantom{-}I & &\\
 & & &  &  & & & & &  & &\textcolor{green}{I} & \textcolor{blue}{\tilde{Z}_3^\top} & & & & & & & &\\
 & & & & & & & & & \phantom{-}\textcolor{orange}{V_3}& \textcolor{red}{N_3}&\textcolor{blue}{\tilde{Z}_3} & & & & & & & &\\
 & & &  &  & I & & & &  -I & & & & & & & & & & &\\
\hline
 & & & & & &  & & & & & & & & \phantom{-}\textcolor{green}{I} & & & \textcolor{orange}{V^\top_4} & -I & & \\
 & & & & & & & &  & & & & & &  &  \textcolor{green}{I}  & & \textcolor{red}{ N^\top_4} & & & \phantom{-}I\\
 & & & & & & & &  &  & & & & &  & &\textcolor{green}{I}  & \textcolor{blue}{\tilde{Z}_4^\top} & & &\\
 & & & & & & & & & &  &  & & & \phantom{-}\textcolor{orange}{V_4} & \textcolor{red}{N_4}& \textcolor{blue}{\tilde{Z}_4}& & &\\
 & & & & & & & &  &  & I& & & & -I & & &\\
\hline
  & & & & & & & & & &  &  & & & & & & & & \phantom{-}\textcolor{green}{I} & -I\\
  & & & & & & & & & &  &  & & & &I & & & & -I
\end{pNiceArray}
\begin{pNiceArray}{c}
\uu_1\\
\zz_1\\
\bla_1\\
\arrayrulecolor[rgb]{0.8,0.8,0.8}
\hline
\vv_1\\
\uu_2\\
\zz_2\\
\bla_2\\
\mmu_1\\
\hline
\vv_2\\
\uu_3\\
\zz_3\\
\bla_3\\
\mmu_2\\
\hline
\vv_3\\
\uu_4\\
\zz_4\\
\bla_4\\
\mmu_3\\
\hline
\vv_4\\
\mmu_4\\
\end{pNiceArray}
=
\begin{pNiceArray}{c}
\\[4em]
\uparrow\\
\\[4em]
\tilde{b}\\
\\[4em]
\downarrow\\
\\[4em]
\end{pNiceArray}
}
\end{align*}
\caption{Block-tridiagonal permutation of the augmented system when $n = 4$.}
\label{matrix:eq:ord:long}
\end{figure}

\subsection{Multigrid}

Having transformed the augmented system with $B=\partial F$ into
the augmented system with $B=\partial\tilde{F}$, and permuted the
transformed system, we now present our preconditioner as geometric multigrid
applied to \eqref{eq:augsysp2-aj}, i.e., Figure~\ref{matrix:eq:ord:long}.

Following \cite{briggs.2000}, we use the letter $I$ for interpolation and
restriction operators.
In contrast with the use of this letter for the identity, we include both
a subscript and a superscript.
The former denotes the resolution of the grid being mapped from and the
latter denotes the resolution of the grid being mapped to.
For instance, $I_{2\Delta t}^{\Delta t}$ is interpolation from time steps
of length $2\Delta t$ to time steps of length $\Delta t$.
When the context is clear, we also use superscripts to distinguish
quantities on different grids.
We use the Greek letter $\Psi$ for smoothers.
Before discussing our multigrid components in detail, we present a standard
two-level multigrid algorithm as Algorithm~\ref{alg:mgrit}.
Its generalization to other cycles is standard.
We recommend \cite{briggs.2000,trottenberg2000multigrid,mccormick1987multigrid}
for more information.

\begin{algorithm}[H]
\caption{Two-Level Multigrid Algorithm for $\tilde{A}y = \tilde{b}$}
\hspace*{\algorithmicindent}  Input: $\tilde A$, $\tilde b$, $y_{\text{in}}$ \\
\hspace*{\algorithmicindent}  Output: $y_{\text{out}}$
\begin{algorithmic}[1]
		\State Pre-Smooth: $y \gets y_{\text{in}}+ \Psi \left(\tilde{b} - \tilde{A} y_{\text{in}}\right)$
		\State Restrict: $r \gets I_{\Delta t}^{2\Delta t} \left(\tilde{b} - \tilde{A} y\right) $
		\State Coarse Solve: $e \gets \left(\tilde{A}^{2\Delta t}\right)^{-1} r  $
		\State Interpolate: $y \gets y + I_{2\Delta t}^{\Delta t} e $
		\State Post-Smooth: $y_{\text{out}} \gets y + \Psi \left(\tilde{b} - \tilde{A} y\right) $
\label{alg:mgrit}
\end{algorithmic}
\end{algorithm}

\subsubsection{Interpolation and Restriction}

We interpolate and restrict all variables except the control using linear interpolation and injection, respectively.
Let $n$ be divisible by two.
Then, e.g., for the state, we have
\begin{itemize}
\item
\emph{interpolation:}$\quad u_{2k-1}^{\Delta t}=u_k^{2\Delta t}\quad\text{and}\quad u_{2k}^{\Delta t}=\frac{1}{2}(u^{2\Delta t}_{k+1}+u_k^{2\Delta t})$
\item
\emph{restriction:} \hspace{1.2em}$\quad u_k^{2\Delta t} = u_{2k-1}^{\Delta t}$
\end{itemize}
for $k = 1,\ldots, \frac{n}{2}$.
Since we take the control to be piecewise constant over each time step, we define
\begin{itemize}
\item
\emph{interpolation:}$\quad z_{2k}^{\Delta t}=z_{2k-1}^{\Delta t}=z_{k}^{2\Delta t}$
\item
\emph{restriction:} \hspace{1.3em}$z_k^{2\Delta t} = \frac{1}{2}(z_{2k}^{\Delta t}+z_{2k-1}^{\Delta t})$
\end{itemize}
for $k=1,\ldots,\frac{n}{2}$.
Like our choice of time stepping scheme \eqref{eq:theta-aj}, these choices of interpolation and restriction operators are in the interest of concreteness.
Our parallel-in-time multigrid preconditioner supports alternatives.

\subsubsection{Smoother}
\label{sub:sec:smootheners}

We denote the strictly lower block triangular, block diagonal, and strictly upper block triangular parts of $\tilde{A}$ as
$$L_{\tilde{A}},\quad D_{\tilde{A}}, \quad\text{and}\quad U_{\tilde{A}} = L_{\tilde{A}}^\top,$$
respectively.
Our parallel-in-time multigrid preconditioner uses the block Jacobi iteration as its smoother $\Psi$.
For the linear system \eqref{eq:augsysp2-aj}, the iteration is
\begin{align}\label{eq:jacobi-aj}
\begin{split}
  y^{\ell} &= y^{\ell-1}+D_{\tilde{A}}^{-1}\left(\tilde{b}- \tilde{A}y^{\ell-1}\right) \\
             &= D_{\tilde{A}}^{-1}\left[\tilde{b}- \left(L_{\tilde{A}}+L_{\tilde{A}}^\top\right)y^{\ell-1}\right],
\quad \ell = 1,2,\ldots
\end{split}
\end{align}
The structure of $\tilde{A}$ is such that \eqref{eq:jacobi-aj} splits $D_{\tilde{A}}$ from the $u_k=v_k$ conditions that couple time domains.
Therefore, the iteration parallelizes trivially in time.
For \eqref{eq:jacobi-aj} to be mathematically valid, $D_{\tilde{A}}$ must be nonsingular.
The following theorem establishes a sufficient condition for this to be the case.
\vspace{-1em}
\begin{theorem}\label{thm:invertibility}
Suppose the $N_k$ matrices in \eqref{eq:blocks-aj} are invertible for all $k=1,2,\ldots,n$.
Then $D_{\tilde{A}}$ is nonsingular.
\end{theorem}
\vspace{-1em}
\begin{proof}
The theorem follows from an augmented system matrix having full rank when $B$ does.
See, e.g., \cite[Theorem~3.1]{benzi.2005}.
Dropping the off-diagonal identity blocks in Figure~\ref{matrix:eq:ord:long} amounts to dropping the identity blocks on the anti-diagonal of Figure~\ref{matrix:eq:long:v}.
Upon doing so, the remainder of the $B$ matrix in Figure~\ref{matrix:eq:long:v} has full rank when the $N_k$ matrices are invertible.
\end{proof}

\subsubsection{Control Scaling}

Numerical experiments suggest a significant benefit to upweighting the
control inner product relative to those for the other variables,
resulting in an effective block Jacobi smoother.
Specifically, having introduced virtual variables, we now have in place
of~\eqref{eq:ipx-aj}
\begin{align*}
&(x,\tilde{x})_X = (u,\tilde{u})_U + (v,\tilde{v})_U + \gamma^2(z,\tilde{z})_Z
\quad\text{for all}\quad x, \tilde{x}\in U\times U\times Z,
\end{align*}
and we find it advantageous to choose $\gamma\gg 1$.

To motivate choosing $\gamma$ in this manner, we find it convenient to
interchange the positions of the controls and adjoints
in~Figure~\ref{matrix:eq:ord:long} as well as the positions of the
constraint Jacobians and the control identities.
This reordering preserves the block tridiagonal structure
of~Figure~\ref{matrix:eq:ord:long} but modifies the $k$th interior
block row of the tridiagonal to have the structure
\begin{align}\label{eq:blockrow}
\begin{bmatrix} L_k & D_k & L_k^\top \end{bmatrix},
\quad
  L_k = \left[\begin{array}{c|c|c}
&\phantom{A} &\phantom{A} \\
\hline
\pi & & \\
\hline
& & \\
\end{array}\right]
\quad
\text{and}
\quad
  D_k = \left[\begin{array}{c|c|c}
I & \alpha^\top & \\
\hline
\alpha & & \frac{1}{\gamma}\beta \\
\hline
& \frac{1}{\gamma}\beta^\top & I
\end{array}\right].
\end{align}
In \eqref{eq:blockrow}, we have abstracted the variables (i.e., the columns)
into three sets: ($i$) states and virtual states, ($ii$) adjoints,
and ($iii$) controls.
The $\alpha$ block is the Jacobian of the constraint with respect to the first
set of variables at time step $k+1$, modulo the inter-time step coupling
block~$\pi$, and the $\beta$ block is the control Jacobian at time step $k+1$.
Observe that sending $\gamma\to\infty$ fully decouples the first and the second
sets of variables (the states, virtual states and the adjoints) from the third
(the controls).
This decoupling is within and across all block rows of the augmented system
matrix, which permits us to permute the decoupled augmented system matrix into
the form
\begin{align}\label{eq:augsys.block}
\begin{bmatrix} J \\ & I \end{bmatrix}
\end{align}
The block $J$ represents the interactions between states, virtual states, and
adjoints, and $I$ is an identity of the appropriate dimension for the controls.

We motivate a large $\gamma$ through this block diagonal decoupling
with $\gamma\to\infty$.
The block Jacobi iteration for \eqref{eq:augsys.block} should be more effective
than the small-$\gamma$ cases where the $(1,2)$ and $(2,1)$
blocks are generally nonzero.
This is because the block Jacobi iteration for \eqref{eq:augsys.block}
amounts to the block Jacobi iteration for the smaller block $J$ instead of
the interdependent update of all variables.
In Appendix~\ref{a:sparse}, we provide further justification of this intuition.

\subsubsection{Coarse-Grid Solver}\label{sec:coarse-grid}

As the coarse grid solver of our preconditioner, we use GMRES
preconditioned by the symmetric block Gauss-Seidel (SGS) iteration.
For the linear system \eqref{eq:augsysp2-aj}, the iteration is
\begin{align*}
  &y^{\ell} = y^{\ell-1}
+ \left[\left(D_{\tilde{A}}+L_{\tilde{A}}\right)D_{\tilde{A}}^{-1}\left(D_{\tilde{A}}+L_{\tilde{A}}^\top\right)\right]^{-1}
\left(\tilde{b}- \tilde{A}y^{\ell-1}\right),
\quad \ell = 1,2,\ldots
\end{align*}
Unlike the block Jacobi iteration, SGS does not parallelize in time,
due to the presence of the $L_{\tilde{A}}$ and $L_{\tilde{A}}^\top$ blocks.
However, as we demonstrate next, SGS is extremely effective in
solving the coarse-grid optimal control problems.

In a preview of our numerical experiments---see Section~\ref{sec:numerics}
for their detailed descriptions---to explain our choice of SGS as the
coarse-grid solver, Table~\ref{tab:fp-gmres} documents the GMRES iterations
needed to solve augmented systems within the SQP optimization
method~\cite{heinkenschloss.2014} for different fixed-point
block iteration schemes \emph{used as preconditioners} (i.e., not as
coarse-grid solvers).
In addition to block Jacobi and SGS, defined earlier, we consider forward
block Gauss-Seidel (FGS) and backward block Gauss-Seidel (BGS) iterations.
\begin{table}[htb!]
\begin{tabular}{c c  cc c cc c cc c cc}
\hline
&& \multicolumn{2}{c}{{Jacobi}} && \multicolumn{2}{c}{{FGS}} && \multicolumn{2}{c}{{BGS}} && \multicolumn{2}{c}{{SGS}} \\
{$n$}  && {SQP} & {LS}  && {SQP} & {LS} && {SQP} & {LS} && {SQP} & {LS} \\
\hline
64   && 8 & {94.39}  && 8 & {56.34}   && 8 & {40.07}   && 9  & {2.66} \\
128  && 8 & {178.99}  && 8 & {110.51}  && 8 & {74.79}   && 8  & {2.72}  \\
256  && 8 & {338.46}  && 8 & {219.45}  && 8 & {141.35}  && 8  & {2.51}  \\
512  && 8 & {645.65}  && 8 & {436.10}  && 8 & {274.39}  && 10 & {2.16}  \\
\hline
\end{tabular}
  \caption{Performance of fixed-point iterations for the van der Pol control example,
  as right preconditioners for GMRES. Here, $n$ is the number of time steps,
  {SQP} is the number of SQP iterations and {LS} (an abbreviation for linear solves)
  is the average number of GMRES iterations per call.
  The control scale factor is set to $\gamma = 100$.
  The relative residual tolerance for GMRES is set to the minimum of~$10^{-6}$
  and the dynamically adjusted tolerances given in~\cite{heinkenschloss.2014}.
  We recognize the remarkable performance advantage of SGS, which is why
  we choose it as our coarse-grid preconditioner.
  }
  \label{tab:fp-gmres}
\end{table}
We note that SGS performs remarkably well as an augmented system
preconditioner.
Its performance is independent of the number of time steps, $n$.
The iteration numbers are less than 3 on average, while the iteration numbers
for block Jacobi, FGS and BGS increase significantly as $n$ increases.
Even though we do not expect the larger numbers of time steps studied in
Table~\ref{tab:fp-gmres} to be used on practical coarse grids, due to its
efficiency and robustness we choose SGS-preconditioned GMRES as our
coarse-grid solver.
The time serialization inherent to SGS is not expected to impact
the parallel scalability of the multigrid solver as long as the coarse grid
is relatively small.
As a final remark, in practice we do not solve the coarse systems to
machine precision, or even to the tolerance from Table~\ref{tab:fp-gmres};
we choose a relative residual tolerance of $10^{-3}$.
Thus, the application of our multigrid preconditioner is nonlinear and
requires special care for use in a Krylov solver.
We resolve this by employing flexible GMRES~\cite{saad.1993} in the outer loop.

%
%
%
\section{Numerical Results}
\label{sec:numerics}

We demonstrate the utility of our parallel-in-time multigrid preconditioner
in the context of augmented systems arising in a nonlinear equation solver
and a nonlinear optimization algorithm.
Specifically, we solve
\begin{itemize}
  \item an optimal control problem with the van der Pol oscillator equation;
  \item an optimal control problem with the viscous Burgers' equation; and
  \item the inviscid Burgers' equation (nonlinear equation). 
\end{itemize}
To solve the nonlinear control problems, we use the matrix-free
SQP algorithm~\cite{heinkenschloss.2014}.
This algorithm is a composite-step method, involving the computation of
a so-called \emph{quasinormal step} that solves the nonlinear constraint
equation through a sequence of constraint linearizations, and a
\emph{tangential step} that improves optimality while staying close to the
nullspace of the linearized constraint equation.
Both steps require approximate solutions of augmented systems; in the
form of the minimum-norm problem~\eqref{eq:mn-aj} for the quasinormal step,
and the orthogonal projection~\eqref{eq:ls-aj} to maintain the proximity of the
tangential step to the linearized constraint nullspace.
This orthogonal projection is used repeatedly within a projected conjugate
gradient (CG) subalgorithm, which dominates the computational cost of the
SQP method.
To solve the nonlinear equation (and, similarly, to compute the quasinormal
step in the SQP method), we combine Powell's dogleg
method~\cite{powell1970hybrid,powell1970new} with augmented system solves
to generate its Gauss-Newton steps.

We solve augmented systems using flexible GMRES (FGMRES).
The flexible iteration is needed because our multigrid preconditioner is
nonlinear, due to the inner GMRES-based SGS-preconditioned
coarse-grid solve with a relative tolerance of $10^{-3}$.
Within the SQP algorithm, we terminate the outer FGMRES iterations using a
relative residual stopping tolerance set to the minimum of~$10^{-6}$ and the
linear solver stopping conditions (LSSCs)
given in~\cite[LSSC 4.3, 4.14, 4.15, 4.17, and 4.22]{heinkenschloss.2014}.
For our dogleg nonlinear solver, the stopping condition LSSC 4.3
from~\cite{heinkenschloss.2014} applies directly.
We note that the LSSCs from~\cite{heinkenschloss.2014} are based on the
progress of the optimization algorithm---resulting in nonuniform relative
residuals that are typically significantly larger than~$10^{-6}$---however,
here, we are primarily interested in studying the performance of the linear
solver, which is why we force its relative stopping tolerances to a uniform
value of~$10^{-6}$.
In almost all augmented system solves, this uniform value is the effective
stopping tolerance.

We focus our numerical studies on the algorithmic performance
of the standard V multigrid cycle, with a coarsening factor of two.
We experimented with W and F cycles, including larger coarsening factors,
however we found no significant advantages to them in our numerical examples.
A fundamental challenge is that the multigrid cycles themselves
serialize computations and that the more complex cycles may limit the
potential for time parallelization in high-performance implementations. 
For instance, while a level-4 V cycle takes 6 sequential steps, level-4
F and W cycles take 18 and 28 sequential steps, respectively.

\subsection{Optimal Control of the Van der Pol Oscillator} \label{sec:numerics-vanderpol}

We examine the control problem
\begin{subequations}
\label{prob:ode}
\begin{align}
\mathop{\text{minimize}}_{u,z}
  &\;\; \frac{1}{2} \int_0^T (u_1(t) - u_{1,d}(t))^2  + (u_2(t) - u_{2,d}(t))^2  dt \nonumber \\
  &\;\; \qquad + \frac{\alpha}{2}\int_0^T z_1^2(t) + z_2^2(t)  dt \\
\text{subject to}
&\;\; \frac{d u_1(t)}{dt} = u_2(t) + z_1(t), \quad t \in (0,T), \\
&\;\; u_1(0) = u_{1,\text{init}},\\
&\;\; \frac{d u_2(t)}{dt} = \mu( 1 - u_1^2(t))  u_2(t)  - u_1(t)+ z_2(t), \quad t \in (0,T),\\
&\;\; u_2(0) = u_{2,\text{init}},
\end{align}
\end{subequations}
where $u_1(t)$, $u_2(t)$ are the state variables,
$z_1(t)$, $z_2(t)$ are the control variables,
and $u_{1,d}(t)$, $u_{2,d}(t)$ are the target data.
The constraint can be viewed as a two-dimensional version of
the van der Pol oscillator equation with a damping
parameter~$\mu$~\cite{james1974time,guckenheimer1980dynamics}.
In our numerical experiments we use
$ u_{1,\text{init}} =  u_{2,\text{init}} = 1$ and generate the
data $u_{1,d}(t)$ and $u_{2,d}(t)$ by
propagating the ODE with $\mu = 0$ and $T = 8$,
creating a circular orbit. 
We then change the damping to $\mu = 0.1$, which generally forces
a non-circular trajectory.
Our goal is to control the trajectory to be more circular,
hence the choice of our data. 
Similar examples have been studied in the control
literature~\cite{diaz2009nonlinear,james1974time,chagas2012optimal}.
We set the penalty parameter to $\alpha =  0.1$ and
use the trapezoidal time stepping rule, i.e., $\theta = 0.5$.

\begin{figure}[htb!]
\centering
\includegraphics[width=0.6\textwidth]{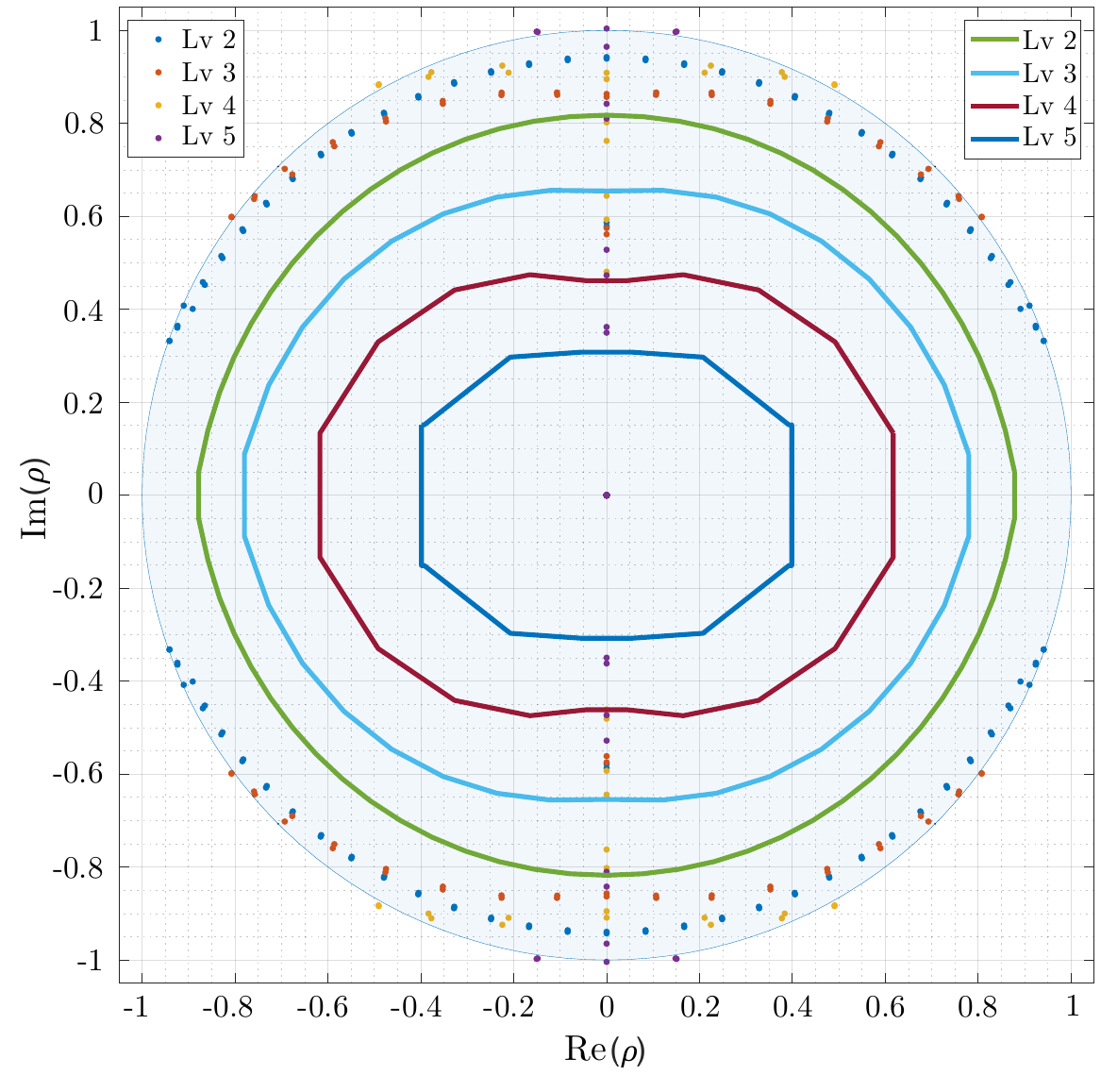}
\caption{Eigenvalues $\rho$ of the block Jacobi iteration matrix,
  $I - D_{\tilde{A}}^{-1}\tilde{A}$,
  at different grid levels, for the 64 time-step discretization of the
  van der Pol control example at the first SQP iteration.
  Different colors indicate different grid levels from
  Table~\ref{tab:spectralradii}.
  \emph{Dots}: Eigenvalues with $\gamma=1$.
  \emph{Bold curves}: Outer boundaries of the eigenvalue regions for $\gamma=100$.
  }
\label{fig:spectrum}
\end{figure}

We highlight the importance of the control scale factor~$\gamma$
for the effectiveness of the block Jacobi smoother.
Table~\ref{tab:spectralradii} documents the spectral radii of the
block Jacobi iteration matrix
$I - D_{\tilde{A}}^{-1}\tilde{A}$
for $\gamma=1$ and $\gamma=100$.
We observe that for $\gamma=1$, spectral radii are either
slightly below or slightly above~1.
In contrast, for $\gamma=100$, spectral radii are smaller
than~1 and decrease with increasing grid levels, tending to
to approximately 0.4 (in bold) at the final grid level.
\begin{table}[htb!]
\begin{tabular}{c c cc c cc c cc}
\hline
&& \multicolumn{2}{c}{{$n=64$}} && \multicolumn{2}{c}{{$n=128$}} && \multicolumn{2}{c}{{$n=256$}} \\
{Lv} && {$\gamma=1$} & {$\gamma=100$} && {$\gamma=1$} & {$\gamma=100$} && {$\gamma=1$} & {$\gamma=100$}\\
\hline
2  && 0.9980 & 0.8797     && 0.9990 & 0.9358     && 0.9995 & 0.9668 \\
3  && 1.0055 & 0.7847     && 1.0037 & 0.8811     && 1.0020 & 0.9374 \\
4  && 1.0109 & 0.6315     && 1.0109 & 0.7831     && 1.0067 & 0.8818 \\
5  && 1.0088 & \bf 0.4264 && 1.0164 & 0.6235     && 1.0140 & 0.7820 \\
6  && --    & --          && 1.0153 & \bf 0.4069 && 1.0192 & 0.6185 \\
7  && --    & --          && --     & --         && 1.0200 & \bf 0.3940 \\
\hline
\end{tabular}
  \caption{Spectral radii of the block Jacobi iteration matrix,
  $I - D_{\tilde{A}}^{-1}\tilde{A}$,
  at different grid levels, for the 64, 128 and 256 time-step discretizations
  of the van der Pol control example at the first SQP iteration.
  Here, $n$ is the number of time steps, Lv is the multigrid level,
  and $\gamma$ is the control scale factor.}
  \label{tab:spectralradii}
\end{table}
For the case $n=64$, the eigenvalues of the iteration matrix are
displayed in Figure~\ref{fig:spectrum}.
We note that, for $\gamma=1$, most eigenvalues are clustered around the
boundary (circle) of the unit disk (shaded region).
For $\gamma=100$, the eigenvalue regions are well inside the unit disk
and they shrink with increasing grid levels.

The dramatic effect of a large $\gamma$ is evident from
Table~\ref{tab:vdp-mgit}.
Comparing $\gamma=1$ and $\gamma=100$, the latter significantly reduces the
average number of FGMRES iterations, by up to 24~times, while only doubling
the numbers of SQP and CG iterations.
With $\gamma=100$, the average numbers of FGMRES iterations grow very slowly,
as the logarithm of the numbers of time steps; also see
Figure~\ref{fig:mgit-complexity}.
Overall, our multigrid-in-time preconditioner exhibits excellent algorithmic
scalability with small coarse grids of fixed size.
\begin{table}[htb!]
\begin{tabular}{c c  c  c c c  c  c c c}
\hline
&&& \multicolumn{3}{c}{{$\gamma=1$}} && \multicolumn{3}{c}{{$\gamma=100$}} \\
{$n$}  & {Lv} && {SQP} & {CG} & {LS} && {SQP} & {CG} & {LS} \\
\hline
64   & 4  && 5 & 24 & {27.59}  && 9  & 54 & {4.21}    \\
128  & 5  && 5 & 24 & {55.09}  && 10 & 60 & {5.16}    \\
256  & 6  && 5 & 24 & {80.39}  && 11 & 67 & {5.78}    \\
512  & 7  && 5 & 23 & {101.71} && 11 & 67 & {6.67}    \\
1024 & 8  && 5 & 24 & {175.59} && 11 & 61 & {7.36}    \\
2048 & 9  && 5 & 23 & {176.64} && 10 & 56 & {8.33}    \\
\hline
\end{tabular}
  \caption{Multigrid-in-time performance for the van der Pol control example,
           using the V cycle with four pre-smoothing and four
           post-smoothing steps.  Here, {$n$} is the number of time steps,
           {Lv} is the last multigrid level, {SQP} is the number of
           SQP iterations, {CG} is the number of projected CG iterations, and
           {LS} is the average number of FGMRES iterations per call to FGMRES.}
  \label{tab:vdp-mgit}
\end{table}

\subsection{Solution and Control of Burgers' Equation} \label{sec:numerics-burgers}

We study an optimal control problem involving the viscous Burgers'
equation~\cite{heinkenschloss2008numerical,antoulas2019loewner}.
The continuous optimal control problem is
\begin{subequations}\label{prob:pde}
\begin{align}
\mathop{\text{minimize}}_{u,z}  &\;\; \frac{1}{2} \int_0^T \int_0^1 (u(x,t) - u_{d}(x,t))^2 \ dx dt
                  + \frac{\alpha}{2}\int_0^T \int_0^1 z^2(x,t) \ dx dt\\
\text{subject to}
&\;\; \frac{\partial u(x,t)}{ \partial  t} - \nu \frac{\partial^2  u(x,t)}{ \partial  x^2}
      + \frac{\partial u(x,t)}{ \partial x} u(x,t) = z(x,t) , \nonumber \\
&\;\; \qquad\qquad\qquad\qquad\qquad\quad\;\; (x,t) \in (0,1)\times(0,T), \label{eq:burgers-begin} \\
&\;\; u(0,t) = 0 , \quad t \in (0,T), \\
&\;\; u(1,t) = 0 , \quad t \in (0,T), \\
&\;\; u(x,0) = u_0(x) , \quad x \in (0,1),
\end{align}
with
\begin{align}
u_0(x) =
 \begin{cases}
          1 \quad &\text{if} \, x \in (0,0.5] \\
          0 \quad &\text{if} \, x \in (0.5,1) \\
 \end{cases},
     \quad
     u_d(x,t) = u_0(x),
     \quad t \in (0,T).
\label{eq:burgers-end}
\end{align}
\end{subequations}
In this formulation the functions $u$ and $z$ are the state and control,
respectively, where the control is intended to maintain the initial
condition over the time interval $(0,T)$ for $T=1$.
As in~\cite{heinkenschloss2008numerical}, we choose $\alpha=0.1$.
The state equation is discretized in space using the upwind finite
difference scheme with 512 intervals, over the $[0, 1]$ domain, and
the backward Euler method in time, i.e., $\theta = 1$.
We choose the backward Euler method due to the known issues with the
L stability of the trapezoidal rule (i.e., Crank Nicolson method) in the
context of Burgers' equation~\cite[p. 406--407]{leveque2002finite}. 
For our optimal control studies with the SQP method we set the
viscosity parameter to $\nu = 0.01$.
Additionally, we apply the dogleg nonlinear solver to the inviscid
Burgers' equation, i.e., we
solve~\eqref{eq:burgers-begin}-\eqref{eq:burgers-end} with $\nu = 0$.

Table~\ref{tab:burgers-mgit} documents the performance of our
multigrid-in-time preconditioner with control scalings $\gamma=1$
and $\gamma=100$.
As in Section~\ref{sec:numerics-vanderpol}, $\gamma=100$ results
in significantly better algorithmic performance of our preconditioner,
which exhibits only modest iteration growth as $n$ increases and the
size of the coarse grid is kept small and fixed.
This is likely due to the improved performance of the block Jacobi
smoother.
Additionally, we observe improved performance of the SGS-preconditioned
GMRES solver for the coarse grid, which requires only a single GMRES
iteration for larger numbers of time steps---a six-fold reduction.
Finally, we note that the tradeoff for using $\gamma=100$ is a very
small increase in the number of projected CG iterations within the
SQP algorithm.

\begin{table}[ht]
\begin{tabular}{c c  c  c c c c  c  c c c c}
\hline
&&& \multicolumn{4}{c}{{$\gamma=1$}} && \multicolumn{4}{c}{{$\gamma=100$}} \\
{$n$}  & {Lv} && {SQP} & {CG} & {SGS} & {LS} && {SQP} & {CG} & {SGS} & {LS} \\
\hline
64   & 4  && 7 & 27 & 6.94 & {9.39}  && 6 & 30 & 2.00 & {3.58}    \\
128  & 5  && 7 & 26 & 6.74 & {12.93} && 7 & 34 & 1.54 & {4.75}    \\
256  & 6  && 8 & 32 & 6.25 & {18.94} && 8 & 37 & 1.18 & {6.02}    \\
512  & 7  && 9 & 34 & 6.08 & {27.44} && 8 & 38 & 1.04 & {7.62}    \\
1024 & 8  && 9 & 36 & 6.01 & {39.93} && 9 & 44 & 1.03 & {9.92}    \\
2048 & 9  && 8 & 33 & 5.97 & {57.46} && 8 & 39 & 1.01 & {12.26}   \\
\hline
\end{tabular}
  \caption{Multigrid-in-time performance for the Burgers' control example,
           using the V cycle with four pre-smoothing and four
           post-smoothing steps.  Here, {$n$} is the number of time steps,
           {Lv} is the last multigrid level, {SQP} is the number of
           SQP iterations, {CG} is the number of projected CG iterations
           in the tangential step computation, {SGS} is the average number
           of symmetric block Gauss-Seidel iterations per coarse-grid solve,
           and {LS} is the average number of FGMRES iterations per call.
           The control scale factor is denoted by $\gamma$.}
  \label{tab:burgers-mgit}
\end{table}

\begin{figure}
\begin{center}
  \includegraphics[width=0.9\textwidth]{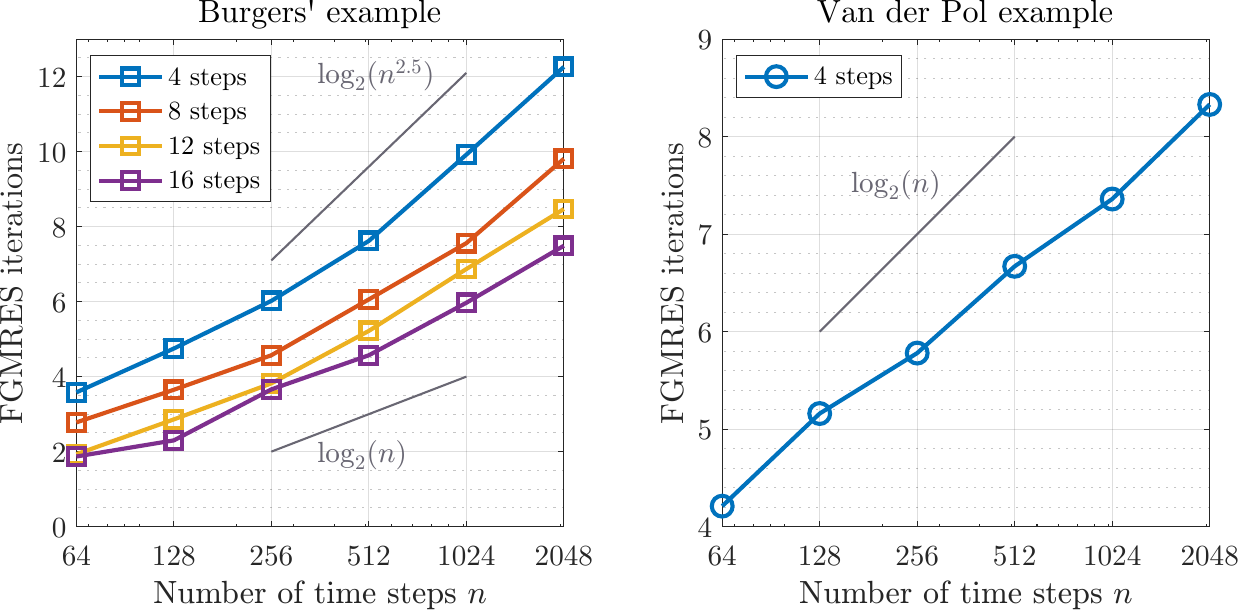}
\end{center}
\caption{Summary of the multigrid-in-time preconditioner performance,
         as illustrated by the average numbers of FGMRES iterations per call.
         \emph{Left}: Burgers' example with increasing numbers of pre-smoothing
         and post-smoothing steps (from Table~\ref{tab:burgers-smoothing}).
         With increased smoothing, iteration numbers tend to~$\log_2(n)$.
         \emph{Right}: For comparison, we include the van der Pol example with
         four pre-smoothing and post-smoothing steps (from Table~\ref{tab:vdp-mgit}).
         Iteration numbers are slightly under~$\log_2(n)$.}
\label{fig:mgit-complexity}
\end{figure}

In Table~\ref{tab:burgers-smoothing} we study the performance of our
preconditioner with increasing numbers of pre-smoothing and post-smoothing
steps of block Jacobi.
Increasing numbers of smoothing steps lead to near $\log_2(n)$ FGMRES
iteration numbers; also see Figure~\ref{fig:mgit-complexity}.
While we do not recommend using large numbers of smoothing steps in
practice (due to their inherent serialization), these performance
numbers confirm conventional multigrid complexity results.
\begin{table}[ht]
\begin{tabular}{c c  c c c c c c}
\hline
&& \multicolumn{6}{c}{Average numbers of FGMRES iterations per call} \\
\begin{tabular}{@{}c@{}}Smoothing \\ steps\end{tabular}  &&
{$n=64$} & {$n=128$} & {$n=256$} & {$n=512$} & {$n=1024$} & {$n=2048$} \\
\hline
4   &&  3.58  &  4.75  &  6.02  &  7.62  &  9.92  &  12.26 \\
8   &&  2.78  &  3.65  &  4.57  &  6.06  &  7.56  &  9.81  \\
12  &&  1.93  &  2.86  &  3.82  &  5.22  &  6.87  &  8.47  \\
16  &&  1.87  &  2.30  &  3.66  &  4.57  &  5.97  &  7.49  \\
\hline
\end{tabular}
\caption{Multigrid-in-time performance for the Burgers' control example,
           using the V cycle with increasing numbers of pre-smoothing and
           post-smoothing steps.
           Here, {$n$} is the number of time steps,
           and `Smoothing steps' denotes the number of pre/post smoothing steps
           in the multigrid algorithm.}
\label{tab:burgers-smoothing}
\end{table}

In our final experiment we turn our attention to the solution of
the inviscid (hyperbolic) Burgers' equation using the dogleg nonlinear
solver with augmented system solves.
We use a simple continuation strategy where the solution for the
viscous case $\nu=10^{-1}$ is used as the initial guess for the viscous
equation with $\nu=10^{-2}$, and so on, until the solution
for $\nu=10^{-4}$ is used as the initial guess for the inviscid case. 
The time grid is fixed at 512 time steps, with seven multigrid levels.
Since there are no concerns with the performance of the nonlinear scheme
with large $\gamma$, in contrast to the solution of the optimal control
problem, we choose $\gamma=10^{5}$.
In Table~\ref{tab:hyperbolic}, we observe a single SGS-preconditioned
GMRES iteration on the coarse grid, for all viscosity parameters.
Furthermore, we observe a moderate increase of outer FGMRES iterations
until $\nu=10^{-4}$, at which point the FGMRES iterations decrease.
Following this, solving the inviscid equation takes an average of 12.8
FGMRES iterations.
This is an interesting result, as the hyperbolicity of the partial
differential equation does not significantly hinder the performance of
our multigrid preconditioner, including the coarse-grid solver, further
reinforcing the case for multigrid in time for augmented systems.

\begin{table}[ht]
\begin{minipage}{0.5\textwidth}
\begin{tabular}{c  c   c c c}
\hline
{$\nu$}  && {Dogleg} & {SGS} & {LS} \\
\hline
$10^{-1}$  && 2 & 1.00 & {2.00}     \\
$10^{-2}$  && 7 & 1.00 & {5.62}     \\
$10^{-3}$  && 9 & 1.00 & {11.00}    \\
$10^{-4}$  && 6 & 1.00 & {7.33}    \\
0          && 5 & 1.00 & {12.80}   \\
\hline
\end{tabular}
\end{minipage}
\begin{minipage}{0.5\textwidth}
\hspace{3em}
\includegraphics[width=0.8\textwidth]{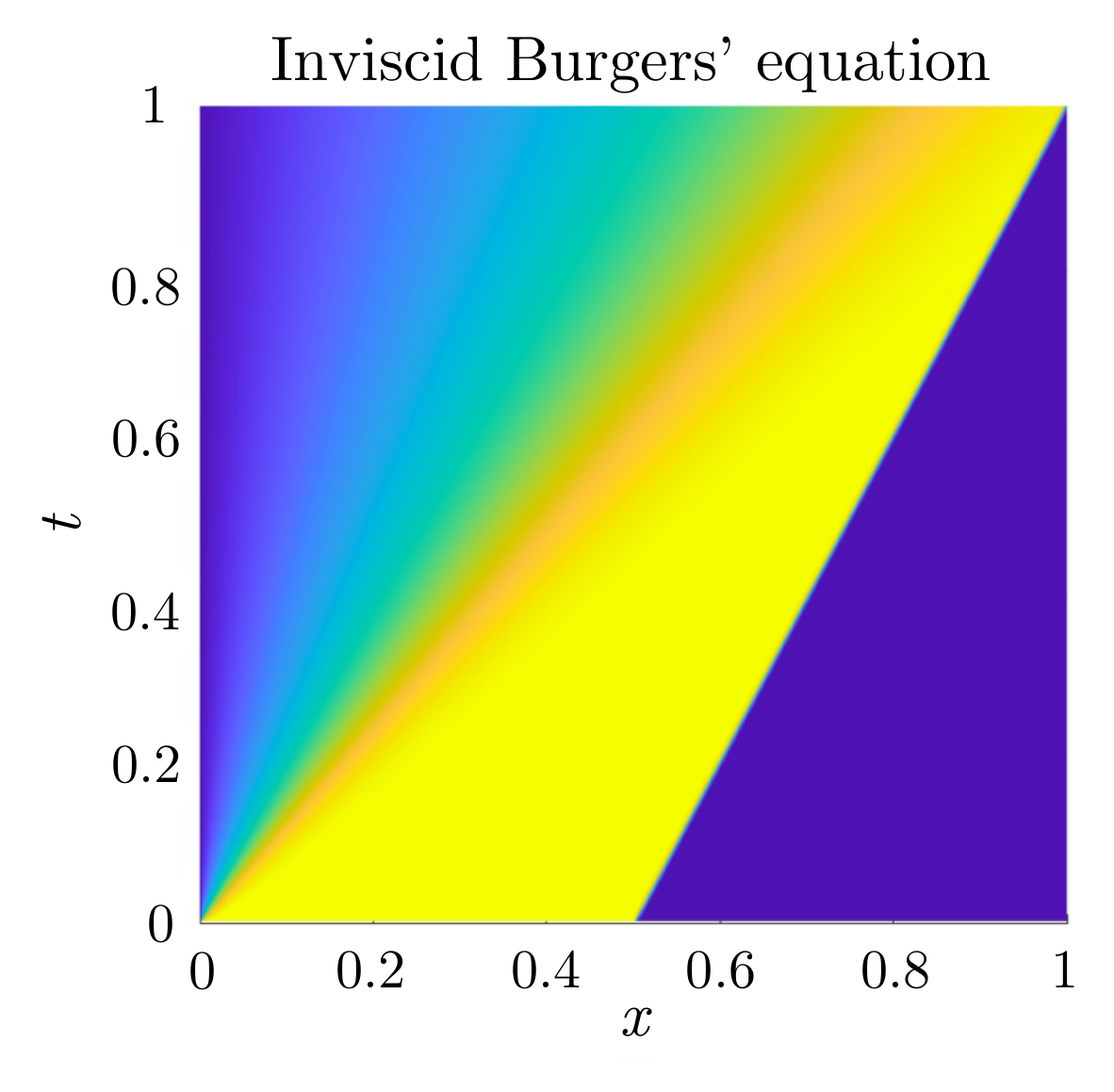}
\end{minipage}
  \caption{\emph{Left:} Multigrid-in-time performance for Burgers' nonlinear equation,
           using the V cycle with four pre-smoothing and four
           post-smoothing steps, as the viscosity decreases from $\nu=10^{-1}$
           to $\nu=0$.  The time grid is fixed at 512 time steps, and the
           number of levels is 7.  Here, {Dogleg} is the number of the
           dogleg nonlinear solver iterations, {SGS} is the average number
           of symmetric block Gauss-Seidel iterations per coarse-grid solve,
           and {LS} is the average number of FGMRES iterations per call.
           \emph{Right:} Solution of the inviscid Burgers' equation with
           the step-function initial condition.}
  \label{tab:hyperbolic}
\end{table}

\subsection{Parallelizability}

We close this section with a few notes on our preconditioner's potential for
parallelization in the optimal control setting.
As a performance baseline for serial execution, we take SGS-preconditioned GMRES,
which exhibited scalable behavior on the van der Pol example, with approximately
three iterations needed for convergence to a relative tolerance of $10^{-6}$.
We observed similar behavior on Burgers' optimal control example.
As the base ``serial units'' of computation, we consider augmented system
solves that are executed sequentially on the time subdomains.
An SGS iteration amounts to $2n$ such solves, due to its forward and backward
block Gauss-Seidel structure.
Let $\mathcal{N}_{\text{LS}}^{\text{ser}}$ denote the number of
SGS-preconditioned GMRES iterations.
Thus, the number of serial units of computation for our serial baseline is
$$2 \mathcal{N}_{\text{LS}}^{\text{ser}} n.$$

To analyze the performance of our augmented system solver, we use the following
notation and assumptions.
\begin{itemize}
  \item Let $\mathcal{N}_{\text{LS}}^{\text{par}}$ denote the number of
        FGMRES iterations, $\mathcal{N}_{\text{Lv}}$ denote the number of
        V multigrid levels, and $\mathcal{N}_{\text{pre}}$ and
        $\mathcal{N}_{\text{post}}$ denote the numbers of pre-smoothing and
        post-smoothing steps with the block Jacobi iteration.
  \item We assume unlimited parallel resources.  Under this assumption,
        the cost of each smoothing step is one serial unit of computation,
        because all time-subdomain augmented systems, within a single smoothing
        step, can be solved at the same time.
  \item We assume negligible cost of the coarse-grid solve, which is reasonable
        for sufficiently small coarse grids, i.e., many multigrid levels and/or
        large coarsening factors.
        We also neglect the cost of interpolation and restriction, as they are
        trivial compared to the solutions of linear systems.
        Similarly, we neglect the cost of applying the augmented system
        operators.
\end{itemize}
Sequential operations in our solver comprise the FGMRES iterations, which
encompass the sequential traversal of the multigrid levels (coarsening followed
by refinement), each of which includes pre-smoothing and post-smoothing steps,
executed one after another.
Therefore, under the above assumptions, the number of serial units of
computation for our multigrid-based solver is
$$
  \mathcal{N}_{\text{LS}}^{\text{par}}
  2(\mathcal{N}_{\text{Lv}}-1)
  (\mathcal{N}_{\text{pre}} + \mathcal{N}_{\text{post}}) ,
$$
resulting in a potential parallel speedup of
$$
  \frac{\mathcal{N}_{\text{LS}}^{\text{ser}}}
  {\mathcal{N}_{\text{LS}}^{\text{par}}
  (\mathcal{N}_{\text{Lv}}-1)
  (\mathcal{N}_{\text{pre}} + \mathcal{N}_{\text{post}})}
  \cdot n .
$$
As an example, consider the observed (approximate) and selected values of
$\mathcal{N}_{\text{LS}}^{\text{ser}}=3$,
$\mathcal{N}_{\text{LS}}^{\text{par}}=6$,
$\mathcal{N}_{\text{Lv}}=6$, and
$\mathcal{N}_{\text{pre}}=\mathcal{N}_{\text{post}}=4$;
here, the potential parallel speedup is $n/80$.

%
%
%
\section{Conclusion}

Our parallel-in-time multigrid preconditioner for augmented systems
accelerates the numerical solution of optimal control problems involving
the van der Pol oscillator and the viscous Burgers' equation, as well as
the solution of the inviscid Burgers' equation.
Unlike several other parallel-in-time methods for optimal control and
even other preconditioners for augmented systems \cite{rees.2010}, our
preconditioner simultaneously solves for the augmented state and
adjoint, in parallel.
The preconditioner lifts the augmented system to a higher-dimensional
space through virtual interface variables that promote time-domain
decoupling, followed by applying geometric multigrid.
The multigrid consists of a block Jacobi smoother, which parallelizes
trivially in time, and a coarse grid solver involving GMRES preconditioned
by the symmetric block Gauss-Seidel iteration.
Our preconditioner preserves the structure of augmented systems on the
time subdomains, and, among all parallel-in-time methods for optimal
control, is most closely connected with multigrid.
Three avenues of future work are of particular interest to us.
The first is a high-performance implementation, to fully explore
the preconditioner's potential for parallelization.
The second is to better understand the usefulness of our work in solving
hyperbolic problems (Table~\ref{tab:hyperbolic}), including comparisons
with the state of the art~\cite{howse2019parallel}.
The third is to assess the effectiveness of our work in accelerating
general multiple shooting methods \cite{ascher.1995}, with more
sophisticated, higher-order time stepping schemes.

%
%
%
\section*{Acknowledgments}

This paper describes objective technical results and analysis. Any subjective views or opinions that might be expressed in the paper do not necessarily represent the views of the U.S. Department of Energy or the United States Government.

This article has been authored by an employee of National Technology \& Engineering Solutions of Sandia, LLC under Contract No. DE-NA0003525 with the U.S. Department of Energy (DOE). The employee owns all right, title and interest in and to the article and is solely responsible for its contents. The United States Government retains and the publisher, by accepting the article for publication, acknowledges that the United States Government retains a non-exclusive, paid-up, irrevocable, world-wide license to publish or reproduce the published form of this article or allow others to do so, for United States Government purposes. The DOE will provide public access to these results of federally sponsored research in accordance with the DOE Public Access Plan https://www.energy.gov/downloads/doe-public-access-plan .

\begin{appendices}

\section{Control Scaling Motivation}\label{a:sparse}

In this appendix, we motivate a large value of the control scale factor $\gamma$ by showing that the block rows \eqref{eq:blockrow} of the augmented system matrix are no less diagonally dominant in the $\gamma\to\infty$ limit.
Note, however, that while the block Jacobi iteration tends to be more effective the more diagonally dominant the matrix, the spectral radius of the iteration matrix need not be smaller \cite{golub.2013}.
\begin{proposition}
Suppose the assumptions of Theorem~\ref{thm:invertibility} hold with $L_k$ and $D_k$ defined in \eqref{eq:blockrow}.
Then the $\gamma\to\infty$ limit of
\begin{align*}
  \left\|D_k^{-1}\left[\;L_k\;\Big\vert\; L_k^\top\;\right]\right\|_\infty
\end{align*}
is no greater than its finite $\gamma$ counterpart to first order.
\end{proposition}
\vspace{-1em}
\begin{proof}
  Let $\tilde{D}$ be the $\gamma\to\infty$ limit of $D_k$.
A straightforward application of Theorem~\ref{thm:invertibility} establishes that $\tilde{D}$ is invertible.
Let the sparsity pattern of $\tilde{D}^{-1}$ be denoted as $\texttt{sparsity}(\tilde{D}^{-1})$.
Since $\tilde{D}$ is block diagonal,
\begin{align}\label{eq:boost1-aj}
  \texttt{sparsity}(\tilde{D}^{-1}) = \begin{bmatrix}
  \times & \times & \\
  \times & \times & \\
  &  & \times \end{bmatrix}.
\end{align}
Meanwhile, $D_k = \tilde{D} + \delta$, with
\begin{align}\label{eq:boost2-aj}
  \texttt{sparsity}(\delta) = \begin{bmatrix}
  \phantom{\times} \\
  & & \times \\
  & \times & \end{bmatrix}.
\end{align}
The $\delta$ correction is a function of $\gamma$, which we choose to be sufficiently large to ensure $\|\tilde{D}^{-1}\|\|\delta\|\ll 1$.
It follows that
$$D_k^{-1} = \tilde{D}^{-1}(I+\tilde{D}^{-1}\delta)^{-1} \approx \tilde{D}^{-1} - \tilde{D}^{-1}\delta\tilde{D}^{-1},$$
where we take $\approx$ to mean equal to first order in $\delta$.
From \eqref{eq:boost1-aj} and \eqref{eq:boost2-aj}, we have
\begin{align*}
  \texttt{sparsity}(\tilde{D}^{-1}\delta\tilde{D}^{-1}) = \begin{bmatrix}
  \phantom{\times} & & \times \\
  & & \times \\
  \times & \times & \end{bmatrix},
\end{align*}
and matrix multiplication implies
\begin{align*}
  \texttt{sparsity}\left(\;\tilde{D}^{-1}\;\left[\;L_k\;\Big\vert\; L_k^\top\;\right]\right)
= \left[\begin{array}{ccc|ccc}
  \times & \phantom{\times} & \phantom{\times} & \phantom{\times} & \times & \phantom{\times} \\
  \times &                  &                  &                  & \times & \\
         &                  &                  &                  &        & \\
  \end{array}\right]
\end{align*}
with
\begin{align*}
  \texttt{sparsity}\left(\tilde{D}^{-1}\delta\tilde{D}^{-1}\left[\;L_k\;\Big\vert\; L_k^\top\;\right]\right)
= \left[\begin{array}{ccc|ccc}
         & \phantom{\times} & \phantom{\times} & \phantom{\times} & \phantom{\times} & \phantom{\times} \\
         &                  &                  &                  &                  & \\
  \times &                  &                  &                  & \times                  & \\
  \end{array}\right].
\end{align*}
Since the infinity norm is the maximal absolute row sum, we conclude
\begin{align*}
 \left\|\tilde{D}^{-1}\left[\;L_k\;\Big\vert\; L_k^\top\;\right]\right\|_\infty
\le
\left\|\left(\tilde{D}^{-1} - \tilde{D}^{-1}\delta\tilde{D}^{-1}\right)\left[\;L_k\;\Big\vert\; L_k^\top\;\right]\right\|_\infty
\approx 
\left\|D_k^{-1}\left[\;L_k\;\Big\vert\; L_k^\top\;\right]\right\|_\infty,
\end{align*}
as claimed.
\end{proof}

\end{appendices}
%
%
%

%


\bibstyle{sn-mathphys}
\bibliography{ref}


\end{document}